\documentclass[12pt]{amsart}

\usepackage{amsmath,amsthm,amsfonts,amscd,amssymb}
\usepackage[matrix, arrow]{xy}
\usepackage{hyperref}

\title[Homology of free quantum groups]{Homology of free quantum groups}
\thanks{B.C. was supported by an NSERC fellowship, and ANR and JSPS travel funding}
\author {Beno\^\i t Collins}
\address{D\'epartement de Math\'ematique et Statistique, Universit\'e d'Ottawa,
585 King Edward, Ottawa, ON, K1N6N5 Canada
and 
CNRS, Institut Camille Jordan Universit\'e  Lyon 1, 43 Bd du 11 Novembre 1918, 69622 Villeurbanne, France} \email{bcollins@uottawa.ca}
\thanks{J.H. was supported by the DFG Research Traing group 534.}
\author {Johannes H\"artel}
\address{Mathematisches Institut der Universit\"at G\"ottingen, Bunsenstr. 3-5,
37073 G\"ottingen, Germany} \email{johannes@uni-math.gwdg.de}
\thanks{A.T. was supported by funding the DFG Research Traning Groups 534 and 1493, and the CRC "Higher Order Structures in Mathematics", G\"ottingen}
\author {Andreas Thom}
\address{Mathematisches Institut der 
Universit\"at G\"ottingen, Bunsenstr. 3-5,
37073 G\"ottingen, Germany} \email{thom@uni-math.gwdg.de}

\newtheorem{theorem}{Theorem}[section]
\newtheorem{lemma}[theorem]{Lemma}
\newtheorem{e-proposition}[theorem]{Proposition}

\newtheorem{e-definition}[theorem]{Definition\rm}

\newtheorem{theoreme}{Th\'eor\`eme}[section]

\newtheorem{definition}[theoreme]{D\'efinition\rm}

\def\og{\leavevmode\raise.3ex\hbox{$\scriptscriptstyle\langle\!\langle$~}}
\def\fg{\leavevmode\raise.3ex\hbox{~$\!\scriptscriptstyle\,\rangle\!\rangle$}}

\newcommand{\Tor}{{\rm Tor\,}}
\def\Cz{\mathbb{C}}

\def\Nz{\mathbb{N}}

\def\im{{\rm im}}
\def\Ext{{\rm Ext}}

\begin{document}

\begin{abstract}
We  compute the Hochschild homology of the free orthogonal quantum group $A_o(n)$.
We show that it satisfies Poincar\'e duality and should be considered to be a $3$-dimensional object.
We then use recent results of R. Vergnioux to
derive results about the $\ell^2$-homology of $A_o(n)$ and estimates on
the free entropy dimension of its set of generators. In particular, we show
that the $\ell^2$ Betti-numbers of $A_o(n)$ all vanish and that the free entropy dimension is less than $1$.
{\it This is the english version of a paper to appear in
C. R. Acad. Sci. Paris. \\
To cite this article: B. Collins, J. H\"artel, A. Thom, C. R.
Acad. Sci. Paris, DOI: 10.1016/j.crma.2009.01.021}
\end{abstract}

\maketitle

\section{Introduction and main results}
\label{}

Let $A_o(n)$ be the unital $*$-algebra generated by self-adjoint elements 
$u_{ij}$, $1 \leq i,j \leq n$, with relations:
\[ \sum_{j=1}^n u_{ij} u_{kj} = \delta_{ik},\quad \mbox{and} \quad \sum_{j=1}^n u_{ji} u_{jk} = \delta_{ik}.\]
We call this algebra the \emph{free orthogonal quantum group}, associated to the integer $n \in \Nz$. 
Its abelianization is naturally isomorphic to the algebra of polynomial functions on the compact 
Lie group $O(n)$.
It is endowed with the coproduct $\Delta u_{ij}=\sum_ku_{ik}\otimes u_{kj}$, the counit
$\varepsilon(u_{ij})=\delta_{ij}$ and the antipode $S(u_{ij})=u_{ji}$.
The maps $\Delta:A_o(n)\to A_o(n)\otimes A_o(n)$, $\varepsilon: A_o(n)\to\mathbb{C}$ and
$S:A_o(n)\to A_o(n)^{op}$ are algebra morphisms, therefore they are fully determined 
by the above relations. 

We are in the setting of compact quantum matrix group of Woronowicz \cite{MR901157} therefore 
there exists a unique left and right invariant Haar state, which we denote $\mu$. 
Since $S^2=id$, it follows from Woronowicz's theory that
$\mu$ is a finite trace on $A_o(n)$. It has been proved by Woronowicz that this state
is faithful on $A_o(n)$, and it was proved in \cite{VaVe} that the von Neumann algebra obtained
by the GNS completion of $A_o(n)$ in $\mu$, denoted by $A_o(n)''$,
is actually a $II_1$ factor for $n\geq 3$ 
and that this factor is full and solid. 
These properties are shared by the von Neumann algebra of the free group, therefore it is 
an interesting question whether $A_o(n)''$ is isomorphic to an interpolated free group factor. 
Although we are not able to answer this question, we observe that unlike in the
free group case, all $\ell^2$ Betti numbers vanish. This relies on recent of work of 
R.\ Vergnioux \cite{Ve}, who showed that the first $\ell^2$ Betti number of $A_o(n)$ vanishes.
As far as we know, this is the first example of complete computation of 
$\ell^2$ Betti numbers for a quantum group beyond the case of amenable quantum 
groups (cf Kyed, \cite{Ky}).

We also observe that a corollary of Vergnioux's results following
Connes and Shlyakhtenko \cite{MR2180603} 
is that the non-microstates free entropy dimension of $(A_o(n),\mu)$ is less than $1$,
and this follows from the fact
that its zeroeth and first $\ell^2$ Betti numbers vanish. The problem of computing the microstates
free entropy dimension of $(A_o(n),\mu)$ is still open, as we don't know whether the family
$u_{ij}$ has microstates.

By a slight abuse notation, $\delta_{ij}$ will denote either an element of $\mathbb{C}$ or a multiple of the
identity in $A_o(n)$ depending on the situation. We keep this notation as there will be no ambiguity. 
Let $M_n$ be the $n\times n$ complex matrix algebra over $\mathbb{C}$.
If $N$ is a $\mathbb{C}$-vector space we denote by $M_n N$ the algebraic tensor product of $M_n$ and $N$
over $\mathbb{C}$.
Our first main result is
\begin{theorem}
\label{maintheorem}
Let $n \geq 1$ and $A_o(n)$ be the free orthogonal quantum group.
Let $e_{ij}$ be the canonical matrix basis of the $n\times n$ matrices and 
$r_{ik} = 1 \otimes e_{ik} + \sum_j u_{ij} \otimes e_{kj} \in M_nA_o(n)$.
Then the following sequence $C_*$
of $A_o(n)$ modules
\[ \xymatrix{ 0 \ar[r]& A_o(n) \ar[r]^{\phi_1} & M_n A_o(n)  \ar[r]^{\phi_2} & M_n A_o(n) \ar[r]^{\phi_3} &
A_o(n) \ar[r]^\varepsilon& \Cz \ar[r] & 0}\]
with maps $\phi_1,\phi_2,\phi_3$ defined on the respective module generators 
with 
\begin{eqnarray*}
\phi_1(1) = \sum_{i,j}(u_{ij} - \delta_{ij}) \otimes e_{ij}, \quad
\phi_2(e_{ij}) = r_{ij} \quad \mbox{and} \quad
\phi_3(e_{ij}) = u_{ij} - \delta_{ij}
\end{eqnarray*}
yield a free resolution of the co-unit.
\end{theorem}
We describe some details of the proof in section \ref{proof1}.
Our second main result is 
\begin{theorem} 
\label{maintheorembis}
All $\ell^2$ Betti numbers of $A_o(n)$ vanish. In particular, the free entropy dimension of the
generators of $A_o(n)$ is $1$.
\end{theorem}
We need three ingredients to complete the proof of Theorem \ref{maintheorembis}. 
First we need Theorem \ref{maintheorem};
then we need to observe that our resolution has a self-duality property 
(Theorems \ref{poincare1} and \ref{poincare2}), also known as Poincar\'e Duality,
which allows us to use a result of Thom \cite{Th} to identify the Betti numbers pairwise.
This is done is section \ref{proof2}.
And last we need the recent result of Vergnioux \cite{Ve} mentioned above (section \ref{ell2}).

\section{Description of the resolution of $A_o(n)$}
\label{proof1}

The surjectivity of $\varepsilon$ and the injectivity of $\phi_1$ are obvious.
First of all, let us show that $\im(\phi_3) = \ker(\varepsilon)$. 
The inclusion $\im(\phi_3)\subset \ker(\varepsilon)$ is obvious, so we focus on the proof of 
$\ker(\varepsilon) \subset \im(\phi_3)$. 
Let $a= \sum_{j}a_j  \in A_o(n)$ 
be a sum of monomials with $\varepsilon(a) = 0$,
where a monomial is understood as a product of $u_{ij}$'s. 
This implies that $\sum_j \varepsilon(a_j) =0$ and hence
$a = \sum_j (a_j - \varepsilon(a_j))$. Since $\varepsilon(a_j - \varepsilon(a_j)) =0$, 
it is enough to prove that if $a$ is a monomial 
then $a-\varepsilon (a)\cdot 1\in \im(\phi_3)$. 
If $a$ is of length $0$, there is nothing to prove. 
If $a= a'u_{ij}$ for some $1 \leq i,j \leq n$, then $a=a'(u_{ij} - \delta_{ij}) + a'\delta_{ij}$. 
However, $a'$ has smaller length than $a$. This proves, that we can inductively 
reduce the length of any monomial to $0$. Hence 
$\ker(\varepsilon) = \im(\phi_3)$.

We continue to show that $\ker(\phi_3) = \im(\phi_2)$. For this purpose, we
let $F = \Cz\langle u_{ij}\rangle$ and set 
$R_{ik} = \sum_{j=1}^n u_{ij} u_{kj} - \delta_{ik} \in F$ and
$L_{ik} = \sum_{j=1}^n u_{ji} u_{jk} - \delta_{ik} \in F.$
Moreover, we set
$r_{ik} = 1 \otimes e_{ik} + \sum_j u_{ij} \otimes e_{kj} \in M_nA_o(n)$ and
$l_{ik} = 1 \otimes e_{ki} + \sum_j u_{ji} \otimes e_{jk} \in M_nA_o(n).$
We start with the first intermediate lemma:
\begin{lemma} \label{comput1}
For all $j,k$, the following holds true:
$$r_{jk} = \sum_i u_{ji} l_{ik} \quad \mbox{and} \quad l_{jk} = \sum_{i} u_{ij}r_{ik}.$$
\end{lemma}
Indeed, we compute
$$r_{jk }= \sum_{i} u_{ji} \otimes e_{ki} + \sum_{m} \delta_{jm} \otimes e_{mk} 
= \sum_{i} u_{ji} \otimes e_{ki} + \sum_{i,m} u_{ji} u_{mi} \otimes e_{mk}
= \sum_i u_{ji} l_{ik} \in M_nA_o(n).$$
The other computation is similar. This finishes the proof of lemma \ref{comput1}.
We also need the following lemma:
\begin{lemma} \label{inj}
The linear map $\tilde{\phi_3} \colon M_n(F) \to F$, which is defined by 
$$ (a_{ij})_{i,j =1}^n \stackrel{\tilde\phi_3}{\mapsto} \sum_{ij} a_{ij}(u_{ij} - \delta_{ij}),$$ 
is injective.
\end{lemma}
The proof of this lemma is elementary: assume that $ \sum_{ij} a_{ij}(u_{ij} - \delta_{ij})=0$.
By considering the degree in $u_{ij}$ this implies that $a_{ij}=0$ and this holds true for all
$i,j$.

Now we prove that
$\ker(\phi_3) = \im(\phi_2)$.
Let $a=(a_{ij})_{i,j=1}^n$ by any 
element of $M_n A_o(n)$ which maps to zero under $\phi_3$. Take any lift $\tilde{a} \in M_n F$
of $a \in M_n A_o(n)$. Since $a$ maps to $0 \in A_o(n)$, $\tilde{\phi}_3(\tilde{a})$ has to have a
representation
$$\tilde{\phi}_3(\tilde{a}) = 
\sum_{\alpha} f_{\alpha}R_{i_{\alpha}j_{\alpha}}f'_{\alpha} + 
\sum_{\beta} g_{\beta}L_{i_{\beta}j_{\beta}}g'_{\beta},$$
for some $f_{\alpha}, f'_{\alpha}, g_{\beta}, g'_{\beta} \in F$ and $1 \leq i_{\alpha}, j_{\alpha}, 
i_{\beta}, j_{\beta} \leq n$.
We will construct a pre-image of $\tilde{\phi}_3(\tilde{a})$ by constructing a pre-image of 
each of the summands. If $f'_{\alpha}=1$, then $f_{\alpha}r_{i_{\alpha}j_{\alpha}}$
is a pre-image of $f_{\alpha}R_{i_{\alpha}j_{\alpha}}$. However, if 
$f'_{\alpha} = f''_{\alpha}u_{ij}$, then
$$f_{\alpha}R_{i_{\alpha}j_{\alpha}}f'_{\alpha} =
f_{\alpha}R_{i_{\alpha}j_{\alpha}}f''_{\alpha}(u_{ij} - \delta_{ij}) + 
\delta_{ij}f_{\alpha}R_{i_{\alpha}j_{\alpha}}f''_{\alpha}.$$
Hence we reduced the length of $f'_{\alpha}$ by one on the expense of an additional summand.
However, this additional summand has a pre-image and hence the claim follows by induction.
Since the map $\tilde{\phi_3}$ is injective by Lemma \ref{inj}, we actually reconstructed $\tilde{a}$
which is now by construction of the form:
$$\tilde{a} = \sum_{\gamma} h_{\gamma}R_{i_{\gamma}j_{\gamma}} h'_{\gamma} + 
\sum_{\delta} k_{\delta}R_{i_{\delta}j_{\delta}} k'_{\delta} +
\sum_{i,k} p_{ik} r_{ik} + s_{ik}l_{ik},$$
for some $h_{\gamma}, k_{\delta},  p_{ik}, s_{ik} \in F$, $h'_{\gamma}, k'_{\delta} \in M_n(F)$ and
$1 \leq i_{\gamma},j_{\gamma},i_{\delta},j_{\delta} \leq n$.
Taking the class of $\tilde{a}$ in $M_nA_o(n)$ which is just $a$, we see, using Lemma \ref{comput1}, that
$a \in \im(\phi_2)$. Since $a \in \ker(\phi_3)$ was arbitrary it follows that $\ker(\phi_3) \subset \im(\phi_2)$. 
Again, the reverse inclusion is obvious. This implies that $\ker(\phi_3) = \im(\phi_2)$.

It remains to show that $\ker{\phi_2} = \im{\phi_1}$. The proof of this claim is more involved and uses
a theory of Gr\"obner basis for non-commutative rings. The details can be found in the PhD dissertation
of the second named author \cite{hy}.

\section{Poincar\'e Duality and computation of the $\ell^2$ Betti numbers }
\label{proof2}

In this section we first recall the basic notions of Poincar\'e duality for non-commutative rings, which was
developed in \cite{MR1443171}. For a nice survey and some background we refer to \cite{kra}. In our setup, we can show directly that the resolution is self-dual; indeed $\phi_1^t = \phi_3$ and $\phi_2^t(e_{ij}) = l_{ji}$ as one can easily show. This is sufficient to construct an explicit isomorphism $C_{3-*} \cong \hom_{A_o(n)}(C_{*},A_o(n))$, where $C_* \to \Cz$ denotes the resolution from Theorem \ref{maintheorem}.

\begin{definition} Let $(A,\Delta,S,\varepsilon,\eta)$ be a Hopf-algebra over $\Cz$. 
We define 
$$H_*(A,M) = \Tor_*^{A}(\Cz,M), \quad \mbox{and} \quad H^*(A,M) = \Ext_A^*(\Cz,M),$$
where we consider $\Cz$ as the trivial right $A$-module via the co-unit. Here, $\Tor_*^A(A,?)$ denotes the left-derived functor of the relative tensor product over $A$. Similarly, $\Ext^*_A(A,?)$ denotes the right-derived functor of $\hom_A(A,?)$.
\end{definition}

As usual, we set $A^e = A \otimes_{\Cz} A^{\rm op}$. Since $A$ is a Hopf-algebra $(1 \otimes S) \circ \Delta\colon A \to A^e$ provides a natural embedding. Note that $? \otimes_{A} A^e$ is an exact functor from $A$-modules to $A$-bi-modules which
sends the trivial representation to the trivial bi-module and free modules to free 
bi-modules.
Hence, our resolution of the co-unit provides at the same time a free resolution of the bi-module
$A$ by free $A^e$-modules. Hence, it is equivalent to consider generalized 
group (co)homology in the sense above and Hochschild (co)homology of the Hopf algebra. We state our results only in terms of the objects introduced in the previous definition.

\begin{theorem}
\label{poincare1}
Let $n \geq 2$ and $A_o(n)$ be the free orthogonal quantum group.
\begin{enumerate}
\item $$\Ext^*_{A_o(n)}(\Cz,A_o(n)) = \left\{ \begin{array}{ll} \Cz & *=3 \\
0 & * \neq 3 \end{array} \right.$$
\item There is a natural isomorphism 
$$H_*(A_o(n),M) \cong H^{3-*}(A_o(n),M)$$
for any $A_o(n)$-module $M$.
\item $A_o(n)$ is smooth with $\dim A =3$.
\end{enumerate}
\end{theorem}

Assertion $(i)$ follows straight from the observation the resolution $C_* \to \Cz$ is self-dual. Indeed, the cohomology of the same complex (but upside-down) is clearly $\Cz$ in dimension $3$. The proof of assertion $(ii)$ uses the fact that the $A_o(n)$-modules in the resolution $C_* \to \Cz$ are finitely generated free. Hence,
$$\hom_{A_o(n)}(C_*,M)  \cong \hom_{A_o(n)}(C_*,A_o(n)) \otimes_{A_o(n)} M.$$
The left side computes cohomology (by definition) and the right side computes homology since $C_*$ is self-dual, i.e. $\hom_{A_o(n)}(C_*,A_o(n)) \cong C_{3-*}$.

\begin{theorem} 
\label{poincare2}
Let $A_o(n)$ be the free orthogonal quantum group.
\begin{equation}
H_*(A_o(n),\Cz) = \left\{ \begin{array}{ll} 
\Cz & * =0  \\
\Cz^{\oplus \frac{n(n-1)}{2}} & * =1 \\
\Cz^{\oplus \frac{n(n-1)}{2}} & * =2 \\
\Cz & * = 3 \\
0 & * \geq 4 
\end{array} \right.
\end{equation}
\end{theorem}

Recall that $\Tor_*^{A_o(n)}(\Cz,\Cz)$ can be computed by replacing $\Cz$ with a projective resolution of it,
tensoring over $A_o(n)$ with $\Cz$, and computing the homology. The result can be easily read off the resulting complex,
which is the following:
\[ \xymatrix{
0 \ar[r] & \Cz \ar[r]^{\phi_1} & M_n \Cz \ar[r]^{\phi_2} & M_n \Cz \ar[r]^{\phi_3} & \Cz \ar[r] & 0,}\]
with 
\[\phi_2((a_{ij})_{i,j=1}^n) = \sum_{ij} a_{ij}(e_{ik} + e_{ki}) \quad \mbox{and} \quad 
\phi_1 = \phi_3 = 0.\]
Clearly, $\phi_2(A)=A+A^{t}$ in matrix notation, 
so that the kernel consists of skew-symmetric and the image of all symmetric matrices.
This finishes the proof of Theorems \ref{poincare1} and \ref{poincare2}.

\section{$\ell^2$ homology and free entropy dimension}
\label{ell2}

We recall that the $\ell^2$ Betti numbers of $A_o(n)$ are defined as the dimension 
of the $\ell^2$-homology of the quantum group $H_*(A_o(n); \ell^2(A_0(n),\mu))$. 
It was shown in \cite{Th} that the dimension of the cohomology in the same degree agrees
with the $\ell^2$ Betti number. Theorems  \ref{poincare1} and \ref{poincare2} imply that
$\beta^{(2)}_1(A_o(n),\mu)=\beta^{(2)}_2(A_o(n),\mu),$
and
$\beta^{(2)}_k(A_o(n),\mu)=0$ for all $k\geq 3$.
These assertions follow from Poincar\'e duality and $\beta^{(2)}_0(A_o(n),\mu) =0$. The latter claim follows since $L(A_0(n),\tau)$ is diffuse for $n \geq 2$.

It has been shown recently by R. Vergnioux in \cite{Ve} that
$\beta^{(2)}_1(A_o(n))=0$. His methods uses quantum Cayley graphs.
From his result we can deduce that $\beta^{(2)}_k(A_o(n))=0$ for all $k\geq 3$.
It has been proved by Connes and Shlyakhtenko in \cite{MR2180603}
that the microstate free free entropy dimension of any unital $*$-algebra $A$ with a trace $\tau$ 
satisfies $\delta^*(A,\tau)\leq \beta^{(2)}_1(A,\tau)-\beta^{(2)}_0(A,\tau)+1$.
This implies in our case that  $\delta^*(A_o(n),\mu)=1$. 
We do not know whether $A_o(n)$ has microstates and believe that it is an interesting question. 
If it does, then 
$\delta(A_o(n),\mu)=\delta^*(A_o(n),\mu)=1$.

\section*{Acknowledgements}
The authors would like to thank Roland Vergnioux for sending them an early version of his preprint, 
and the referee for useful comments.
This work was initiated by A.T. and B.C. 
at the Free analysis workshop held at Palo Alto in 2006. We thank this institution 
and the organizers of the workshop for providing us a fruitful working ambiance. 
The three authors are grateful to the ANR Galoisint and C. Kassel from Strasbourg for 
organizing a meeting in February 2008 where they had a chance to get together. 
Part of this work was completed during two visits of B.C. in Universit\"at G\"ottingen, one visit of A.T. in Ottawa University in 2007 and one visit of J.H. in  Coimbra University. The respective visitors are pleased to acknowledge warm hospitality
at the respective visited institutions. J.H. would like to thank Ivan Yudin for many illuminatives talks and hints.

\end{document}